\documentstyle[amssymb,12pt]{article}

\begin{document}

\title{Selfsimilar random fractal measure using contraction method in probabilistic
metric spaces}
\author{J.Kolumb\'an \thanks{%
Babes-Bolyai University, Faculty of Mathematics and Computer Science,
Cluj-Napoca, \newline
e-mail: asoos@math.ubbcluj.ro} and A.So\'os \footnotemark[1] }
\date{}
\maketitle

\newtheorem{th}{Theorem}[section] \newtheorem{prop}{Proposition}[section]
\newtheorem{cor}{Corollary}[section] \newtheorem{lem}{Lemma}[section]

\begin{abstract}
We use contraction method in probabilistic metric spaces to prove
existence and uniqueness of selfsimilar random fractal measures.
\end{abstract}



\section{Introduction}

Contraction methods for proving the existence and uniqueness of nonrandom
selfsimilar fractal sets and measures were first applied by Hutchinson \cite
{H81}. Further results and applications to image compression were obtained
by Barnsley and Demko \cite{BD85} and Barnsley \cite{B88}. At the same time
Falconer \cite{F86}, Graf \cite{G87}, and Mauldin and Williams \cite{MW86}
randomized each step in the approximation process to obtain sefsimilar
random fractal sets. Atbeiter \cite{A91} and Olsen \cite{O94} studied
selfsimilar random fractal measures applying nonrandom metrics. More
recently Hutchinson and R\"uschendorf \cite{HR98a,HR98b,HR00} introduced
probability metrics defined by expectation for random measure and
established existence, uniqueness and approximation properties of
selfsimilar random fractal measures. In these works a finite first moment
condition is essential.

In this paper it will shown that, using probabilistic metric spaces
techniques, we can weak the first moment condition for existence and
uniqueness of selfsimilar measures.

The theory of probabilistic metric spaces, introduced in 1942 by K. Menger
\cite{M42}, was developed by numerous authors, as it can be realized upon
consulting the list of references in \cite{GI89}, as well as those in \cite
{SK83}. The study of contraction mappings for probabilistic metric spaces
was initiated by V. M. Sehgal \cite{S69}, and H. Sherwood \cite{Sh69}.

\bigskip

\section{Selfsimilar random fractal measures}

Recently Hutchinson and R\"{u}schendorf \cite{HR98a,HR98b,HR00} gave a
simple proof for the existence and uniqueness of invariant random measures
using the $L^{q}$-metric, $0<q\leq \infty $. The underlying probability
space for the iteration procedure is generated by selecting independent and
identically distributed scaling laws. Let $(X,d)$ be a complete separable
metric space. A {\em scaling law with weights} ${\bf S}$ is a 2N-tuple $%
(p_{1},S_{1},....,p_{n},S_{N}),\;N\geq 1,$ of positive real numbers $p_{i}$
such that $\sum_{i=1}^{N}p_{i}=1$ and of Lipschitz maps $S_{i}:X\rightarrow X
$ with Lipschitz constant $r_{i}=LipS_{i}$ , i$\in $\{1,...,N\}.

Denote $M=M(X)$ the set of finite mass Borel regular measures on X with the
weak topology. If $\mu \in M$, then the measure ${\bf S}\mu $ is defined by
\[
{\bf S}\mu =\sum_{i=1}^{N}p_{i}S_{i}\mu ,
\]
where $S_{i}\mu $ is the usual push forward measure, i.e.
\[
S_{i}\mu (A)=\mu (S_{i}^{-1}(A)),\,\mbox{for}\,A\subseteq X.
\]
We say $\mu $ {\em satisfies the scaling law {\bf S}} or {\em is a
selfsimilar fractal measure } if ${\bf S}\mu =\mu .$

Let $M_{q}$ denote the set of unit mass Borel regular  measures $\mu $ on X
with finite q-th moment. That is,
\[
M_{q}=\{\mu \in M \,|\,\mu (X)=1,\,\int d^{q}(x,a)d\mu (x)<\infty
\}
\]
for some (and hence any) $a\in X$. Note that, if $p\geq q$ then $%
M_{p}\subset M_{q}.$

The {\em minimal metric} $l_q$ on $M_q$ is defined by
\[
l_q(\mu,\nu)=\inf \{(\int d^q (x,y)d\gamma(x,y))^{\frac 1q \wedge 1}\vert \,
\pi_1\gamma=\mu,\, \pi_2\gamma=\nu\}
\]
where $\wedge$ denotes the minimum of the relevant numbers and $\pi_i\gamma $
denotes the i-th marginal of $\gamma$, i.e. projection of the measure $%
\gamma $ on $X\times X$ onto the i-th component.

We have the following properties of $l_q$ (see \cite{R91}):

a) Suppose $\alpha$ is a positive real, $S: X\rightarrow X$ is Lipschitz,
and $\vee$ denotes the maximum of the relevant numbers. Then for $q>0$ and
for measures $\mu,\nu$ we have the following properties:
\[
l_q^{q\vee 1}(\alpha\mu,\alpha\nu)=\alpha l_q^{q\vee 1}(\mu,\nu),
\]
\[
l_q^{q\vee 1}(\mu_1+\mu_2,\nu_1+\nu_2)\leq l_q^{q\vee
1}(\mu_1,\nu_1)+l_q^{q\vee 1}(\mu_2,\nu_2),
\]
\[
l_q(S\mu,S\nu)\leq (Lip S)^{q\wedge 1}l_q(\mu,\nu)
\]


b) $(M_q,l_q)$ is a complete separable metric space and $l_q(\mu_n,\mu)%
\rightarrow 0$ if and only if

(i) $\mu_n\rightarrow \mu$ (weak convergence) and

(ii) $\int d^q(x,a)d\mu_n(x)\rightarrow \int d^q(x,a) d\mu (x)$ (convergence
of q-th moments).

c) If $\delta_a$ is the Dirac measure at $a\in X$, then
\[
l_q(\mu, \mu(X)\delta_a )=(\int d^q(x,a)d\mu(x))^{\frac 1q \wedge 1},
\]
\[
l_q(\delta_a,\delta_b)=d^{1\wedge q}(a,b).
\]

Let ${\bf M}$ denote the set of all random measures $\mu $ with value in M,
i.e. random variables $\mu :\Omega \rightarrow M.$ Let ${\bf M_{q}}$ denote
the space of random measures $\mu :\Omega \rightarrow M_{q}$ with finite
expected q-th moment i.e.
\begin{equation}
{\bf M_{q}}:=\{\mu \in {\bf M}|\,\mu ^{\omega }(X)=1a.s.,\,E_{\omega
}\int_{X}d^{q}(x,a)d\mu ^{\omega }(x)<\infty \}  \label{csillag}
\end{equation}

The notation $E_\omega$ indicate that the expectation is with respect to the
variable $\omega$. It follows from (\ref{csillag}) that $\mu^\omega\in M_q $
a.s. Note that ${\bf M_p\subset M_q}$ if $q\leq p$.Moreover, since $E^{\frac
1q}\vert f\vert ^q\rightarrow exp(E\log \vert f\vert) $ as $q\rightarrow 0,$
\[
{\bf M_0}:=\cup_{q>0}{\bf M_q}=\{\mu\in {\bf M}\vert \mu_\omega(X)=1\, a.s.
\,, E_\omega\int_X\log d(x,a)d\mu^\omega(x)<\infty\}.
\]
For random measures $\mu,\nu\in {\bf M_q},$ define
\[
l_q^*(\mu,\nu):=\left\{
\begin{array}{lc}
E_\omega^{\frac 1q}l_q^q(\mu^\omega,\nu^\omega), & q\geq 1 \\
E_\omega l_q(\mu^\omega,\nu^\omega), & 0<q<1.
\end{array}
\right.
\]

One can check as in \cite{R91}, that $({\bf M_q{}, l_q^*)}$ is a complete
separable metric space. Note that $l_q^*(\mu,\nu)=l_q(\mu,\nu)$ if $\mu$ and
$\nu$ are constant random measures.

Let ${\cal M}$ denote the class of probability distributions on ${\bf M}$.
i.e.
\[
{\cal M}=\{{\cal D}=dist\mu \,|\,\mu \in {\bf M}\}.
\]
Let ${\cal M}_{q}$ be the set of probability distributions of random
measures $\mu \in {\bf M}_{q}.$ If $q\leq p$ then ${\cal M}_{p}\subset {\cal %
M}_{q}$. Let
\[
{\cal M}_{0}:=\cup _{q>0}{\cal M}_{q}.
\]
The minimal metric on ${\cal M}_{q}$ is defined by
\[
l_{q}^{\ast \ast }({\cal D}_{1},{\cal D}_{2})=\inf \{l_{q}^{\ast }(\mu ,\nu
)|\,\mu \stackrel{d}{=}D_{1},\,\nu \stackrel{d}{=}D_{2}\}.
\]
It follows that $({\cal M}_{q},l_{q}^{\ast \ast })$ is a complete separable
metric space with the next properties:

\[
a) \,l_q^{**}(\alpha {\cal D}_1,\alpha {\cal D}_2)=\alpha l_q^{**}({\cal D}%
_1,{\cal D}_2),
\]
\[
b) \, l_q^{**}({\cal D}_1+{\cal D}_2,{\cal D}_3+{\cal D}_4)\leq {l_q^{**}}^q(%
{\cal D}_1,{\cal D}_3)+{l_q^{**}}^q({\cal D}_2,{\cal D}_4)
\]
for ${\cal D}_i\in {\cal M}_q,\, i=1,2,3,4.$

{\em A random scaling law} ${\bf S}=(p_1,S_1,p_2,S_2,...,p_n,S_N)$ is a
random variable whose values are scaling laws, with $\sum_{i=1}^Np_i=1$ a.s.
We write ${\cal S}=dist {\bf S}$ for the probability distribution determined
by ${\bf S}$ and $\stackrel{d}{=}$ for the equality in distribution.

If $\mu$ is a random measure, then the random measure ${\bf S}\mu$ is
defined (up to probability distribution) by
\[
{\bf S}\mu:=\sum_{i=1}^N p_iS_i\mu^{(i)},
\]
where ${\bf S}, \mu^{(1)},...,\mu^{(N)}$ are independent of one another and $%
\mu^{(i)}\stackrel{d}{=}\mu.$ If ${\cal D}=dist \mu$ we define ${\cal SD}%
=dist {\bf S}\mu.$

We say $\mu$ {\em satisfies the scaling law }${\bf S}$, or is a {\em %
selfsimilar random fractal measure, } if
\[
{\bf S}\mu\stackrel{d}{=}\mu, \;\mbox{or equivalently} \; {\cal SD=D}
\]
and ${\cal D}$ is called a selfsimilar random fractal distribution.

To generate random selfsimilaar fractal measure we use the next {\bf %
iterative procedure }(see \cite{HR98a}):

Fix $q>0$.

Beginning with a nonrandom measure $\mu_0\in M_q$
one iteratively applies iid scaling laws with distribution ${\cal S}$ to
obtain a sequence $\mu_n$ of random measures in ${\bf M_q}$ and a
corresponding sequence ${\cal D}_n$ of distributions in ${\cal \ M}_q$, as
follows:

(i) Select a scaling law ${\bf S}$ via the distribution {\bf S} and define.
\[
\mu_1={\bf S}\mu_0=\sum_{i=1}^n p_i S_i\mu_0, i.e.\,\,\,\mu_1(\omega)={\bf S}%
\mu_0=\sum_{i=1}^n p_i(\omega) S_i(\omega)\mu_0,\, {\cal D}_1\stackrel{d}{=}%
\mu_1,
\]

(ii) Select ${\bf S_1,...,S_n}$ via ${\bf S}$ with ${\bf S^i}=(p_1^i,
S_1^i,...,p_N^i,S_N^i), i\in\{1,2,..., N\}$ independent of each other and of
${\bf S}$ and define
\[
\mu_2:={\bf S}^2\mu_0=\sum_{i,j}p_ip_j^i S_i\circ S_j^i \mu_0,\, {\cal D}_2%
\stackrel{d}{=}\mu_0
\]

(iii) Select ${\bf S}^{ij}=(p_1^i, S_1^{ij},...,p_N^i,S_N^{i,j})$ via ${\cal %
S}$, independent of one another and of ${\bf S^1,...,S^N, S}$ and define

\[
\mu_3={\bf S}^3\mu_0=\sum_{i,j,k} p_ip_j^ip_k^{ij}S_i\circ S_j^i\circ
S_k^{ij}\mu_0,\, {\cal D}_3\stackrel{d}{=}\mu_3,
\]
etc.

Thus $\mu _{n+1}=\sum_{i=1}^{N}p_{i}S_{i}\mu _{n}^{(i)}$ where $\mu
_{n}^{(i)}\stackrel{d}{=}\mu _{n}\stackrel{d}{=}{\cal D}_{n},$ ${\bf S=%
\stackrel{d}{=}{\cal S}},$ and the $\mu _{n}^{(i)}$ and ${\bf S}$ are
independent. It follows that ${\cal D}_{n}={\cal SD}_{n-1}={\cal S}^{n}{\cal %
D_{0}}$, where ${\cal D}_{0}$ is the distribution of $\mu _{0}$. If $\mu
_{0}\in M_{q},$ then ${\cal D}_{0}$ is constant.

The underlying probability space for a.s. convergence is defined above (see
\cite{HR00}).

A {\em construction tree} ( or a construction process ) is a map $\omega
:\{1,...,N\}^{\ast }\rightarrow \Gamma $, where $\Gamma $ is the set of
(nonrandom) scaling laws. A construction tree specifies at each node of the
scaling law used to define constructively a recursive sequence of random
measures. Denote the scaling law of $\omega $ at the node $\sigma $ by the
2N-tuple
\[
{\bf S}^{\sigma }(\omega )=\omega (\sigma )=(p_{1}^{\sigma }(\omega
),S_{1}^{\sigma }(\omega ),...,p_{N}^{\sigma }(\omega ),S_{N}^{\sigma
}(\omega ))
\]
where $p_{i}^{\sigma }$ are weights and $S_{i}^{\sigma }$ Lipschitz maps.
The sample space of all construction trees is denoted by $\tilde{\Omega}.$
The underlying probability space $(\tilde{\Omega},\tilde{{\cal K}},\tilde{P})
$ for the iteration procedure is generated by selecting identical
distributed and independent scaling laws $\omega (\sigma )\stackrel{d}{=}%
{\bf S}$ for each $\sigma \in \{1,...,N\}^{\ast }$. 

In \cite{HR98b} it is proved the following theorem:

\begin{th}
Let ${\bf S}=(p_1,S_1,p_2S_2,...,p_n,S_N)$ be a random scaling
law, with $\sum_{i=1}^Np_i=1$ a.s. Assume $\lambda_q:=E_\omega
(\sum_{i=1}^Np_ir_i^q)<1$ and
\begin{equation}
E_\omega (\sum_{i=1}^Np_id^q(S_ia,a)<\infty \,\mbox{for
some}\,q>0, \,\mbox{and for}\,  a\in X.
 \label{momfelt}
\end{equation}
Then

a) the operator ${\bf S: M_q\rightarrow M_q} $ is a contraction
map with respect to $l_q^*.$

b) If $\mu^*$ is the unique fixed point of ${\bf S}$ and
$\mu_0\in M_p$ (or more generally ${\bf M}_q)$, then
$$E_\omega^{\frac 1q}l_q^q(\mu_n,\mu^*)\leq
\frac{\lambda_q^{\frac{k}{q}}}{1-\lambda_q^{\frac
1q}}E_\omega^{\frac{1}{q}}l_q^q(\mu_0,{\bf S}\mu_0)\rightarrow
0,\,q\geq 1$$
$$E_\omega l_q(\mu_n,\mu^*)\leq \frac{\lambda_q^k}{1-\lambda_q}E_\omega l_q(\mu_0,{\bf S}\mu_0)\rightarrow 0,\, 0<q< 1$$
as $k\rightarrow \infty.$ In particular  $\mu_n\rightarrow \mu^*$
a.s. in the sense of weak convergence of measures.

Moreover, up to probability distribution, $\mu^*$ is the unique
unit mass random measure with $E_\omega \int\ln
d(x,a)d_\mu^\omega <\infty$ which satisfies ${\bf S}$.
\end{th}

Using contraction method in probabilistic metric spaces, instead of
condition (\ref{momfelt}) we can give weaker condition for the existence and
uniqueness of invariant measure. More precisely, in Section 4 we will prove
the following

\begin{th}
\label{sajatmert}
 Let ${\bf S}=(p_1,S_1,p_2S_2,...,p_n,S_N)$ be a
random scaling law, which satisfies  $\sum_{i=1}^Np_i=1$ a.s. and
suppose $\lambda_q:=ess sup( \sum_{i=1}^Np_ir_i^q)<1$ for some
$q>0$. If there exist $\alpha\in M_q$ and a positive number
$\gamma$ such that
\begin{equation}
P(\{\omega\in \Omega\vert\, l_q(\alpha(\omega), {\bf
S}\alpha(\omega))\geq t\})\leq\frac{\gamma}{t},\,\mbox{for all}\,
\,t>0, \label{sajatfelt}
\end{equation}
 then there exists $\mu^*$ such that ${\bf S}\mu^*=\mu^* a.s.$
and  exponentially fast.

Moreover up to probability distribution $\mu^*$ is the unique
unit mass random measure
 which satisfies ${\bf S}$.
\end{th}

{\bf Remark:} If condition (\ref{momfelt}) is satisfied then
condition (\ref {sajatfelt})   hold also. To see this, let $a\in
X$ \ \ and $\ \alpha (\omega ):=\delta _{a}$ for all $\omega \in
\Omega .$ We have:
\[
P(\{\omega \in \Omega |\,l_{q}(\delta _{a}(\omega ),{\bf S}\delta
_{a}(\omega ))\geq t\})=
\]
\[
=P(\{\omega \in \Omega |\,l_{q}(\sum_{i=1}^{N}p_{i}\delta _{a}(\omega
),\sum_{i=1}^{N}p_{i}S_{i}\delta _{a}(\omega ))\geq t\})\leq
\]
\[
\leq P(\{\omega \in \Omega |\,\sum_{i=1}^{N}p_{i}l_{q}(\delta _{a}(\omega ),{%
S_{i}}\delta _{a}(\omega ))\geq t\})=
\]
\[
=P(\{\omega \in \Omega |\,\sum_{i=1}^{N}p_{i}d^{q}(S_{i}a,a))\geq t\})\leq
\frac{1}{t}E_{\omega }(\sum_{i=1}^{N}p_{i}d^{q}(S_{i}a,a))=\frac{\gamma }{t}
\]

However, condition (\ref{sajatfelt}) can  be satisfied also if
\[
E_\omega(\sum_{i=1}^Np_id^q(S_ia,a))=\infty \,\mbox{for all}\,q>0.
\]
Let $\Omega=]0,1]$ with the Lebesque measure, let X be the interval $%
[0,\infty[$ and $N=1$. Define $S:X\rightarrow X$ by $S^\omega(x)=\frac{x}{2}%
+e^{\frac{1}{\omega}}$. This map is a contraction with ratio $\frac 12$. For
$q>0$, the expectation $E_\omega d^q(S0,0)=\infty,$ however
\[
P(\{\omega\in \Omega\vert l_q(S0,0)\geq t\})=\frac{1}{t}
\]
for all $t>0.$



\section{Invariant sets in E-spaces}

\subsection{Menger spaces}

\hspace{0.7cm}Let ${\bf R}$ denote the set of real numbers and ${\bf R_{+}}%
:=\{x\in {\bf R}:x\geq 0\}.$ A mapping $F:${\bf R}$\rightarrow \lbrack 0,1]$
is called a {\em distribution function} if it is non-decreasing, left
continuous with $\inf_{t\in {\bf R}} F(t)=0$ and $\sup_{t\in {\bf R}} F(t)=1$
(see \cite{GI89}). By $\Delta $ we shall denote the set of all distribution
functions $F.$ Let $\Delta $ be ordered by the relation ''$\leq $'', i.e. $%
F\leq G$ if and only if $F(t)\leq G(t)$ for all real t. Also $F<G$ if and
only if $F\leq G$ but $F\not=G$. We set ${\Delta ^{+}}:=\{F\in {\Delta }%
:F(0)=0\}.$

Throughout this paper H will denote the Heviside distribution function
defined by
\[
H(x)=\left\{
\begin{array}{cc}
0, & x\leq 0, \\
1, & x>0.
\end{array}
\right.
\]

\smallskip Let $X$ be a nonempty set. For a mapping ${\cal F}:X\times
X\rightarrow {\Delta ^{+}}$ and $x,y\in X$ we shall denote ${\cal F}(x,y)$
by $F_{x,y}$, and the value of $F_{x,y}$ at $t\in {\bf R}$ by $F_{x,y}(t)$,
respectively. The pair $(X,{\cal F})$ is a {\em probabilistic metric space}
(briefly {\em PM space}) if $X$ is a nonempty set and ${\cal F}:X\times
X\rightarrow {\Delta ^{+}}$ is a mapping satisfying the following conditions:

$1^{0}$. $F_{x,y}(t)=F_{y,x}(t)$ for all $x,y\in X$ and $t\in{\bf R;}$

$2^{0}$. $F_{x,y}(t)=1$, for every $t>0$, if and only if $x=y$;

$3^{0}$. if $F_{x,y}(s)=1$ and $F_{y,z}(t)=1$ then $F_{x,z}(s+t)=1.$

\bigskip

A mapping $T:[0,1]\times \lbrack 0,1]\rightarrow \lbrack 0,1]$ is called a
{\em t-norm} if the following conditions are satisfied:

$4^{0}$. $T(a,1)=a $ for every $a\in[0,1];$

$5^{0}$. $T(a,b)=T(b,a)$ for every $a,b\in[0,1]$

$6^{0}$. if $a\geq c $ and $b\geq d$ then $T(a,b)\geq T(c,d);$

$7^{0}$. $T(a,T(b,c))=T(T(a,b),c)$ for every $a,b,c\in[0,1].$

\bigskip

A {\em Menger space} is a triplet $(X,{\cal F},T),$ where $(X,{\cal F})$ is
a probabilistic metric space, where $T$ is a t-norm, and instead of $3^{0}$
we have the stronger condition

$8^{0}$. $F_{x,y}(s+t) \geq T(F_{x,z}(s),F_{z,y}(t))$ for all $x,y,z\in X $
and $s,t\in{\bf R_{+}.}$

\bigskip

The $(t,\epsilon )$-topology in a Menger space was introduced in 1960 by B.
Schweizer and A. Sklar \cite{SK60}. The base for the neighbourhoods of an
element $x\in X$ is given by
\[
\{U_{x}(t,\epsilon )\subseteq X:t>0,\epsilon \in ]0,1[\},
\]
where
\[
U_{x}(t,\epsilon ):=\{y\in X:F_{x,y}(t)>1-\epsilon \}.
\]


\vspace{0.5cm} In 1966, V.M. Sehgal \cite{S69} introduced the notion of a
contraction mapping in PM spaces. The mapping $f:X\rightarrow X$ is said to
be a {\em contraction} if there exists $r\in ]0,1[$ such that
\[
F_{f(x),f(y)}(rt)\geq F_{x,y}(t)
\]
for every $x,y\in X$ and $t\in {\bf R_{+}.}$

\vspace{0.5cm} A sequence $(x_{n})_{n\in {\bf N}}$ from $X$ is said to be
{\em fundamental} if
\[
\lim_{n,m\rightarrow \infty }F_{x_{m},x_{n}}(t)=1
\]
for all $t>0.$ The element $x\in X$ is called {\em limit} of the sequence $%
(x_{n})_{n\in {\bf N}}$, and we write $\lim_{n\rightarrow \infty }x_{n}=x$
or $x_{n}\rightarrow x$, if $\lim_{n\rightarrow \infty }F_{x,x_{n}}(t)=1$
for all $t>0.$ A probabilistic metric (Menger) space is said to be {\em %
complete} if every fundamental sequence in that space is convergent.

Let $A$ and $B$ nonempty subsets of $X.$ The {\em probabilistic
Hausdorff-Pompeiu distance} between $A$ and $B$ is the function $F_{A,B}:%
{\bf R\rightarrow \lbrack 0,1]}$ defined by
\[
F_{A,B}(t):=\sup_{s<t}T(\inf_{x\in A}\sup_{y\in B}F_{x,y}(s),\inf_{y\in
B}\sup_{x\in A}F_{x,y}(s)).
\]

In the following we remember some properties proved in \cite{KS98,KS01}:

\begin{prop}  If $\cal C$ is a nonempty collection of nonempty closed
bounded sets in a Menger space
$(X,{\cal F},T)$ with $T$ continuous, then $({\cal C, F_C}%
,T)$ is also  Menger space, where ${\cal F_C} $ is defined by
${\cal F_C} (A,B):=F_{A,B}$ for all $A,B\in\cal C$ .
\end{prop}

{\bf Proof.} See \cite{KS98,S69}. $\hspace{1cm}\Box$

\vspace{0.5cm}

\begin{prop}
Let $T_m(a,b):=\max\{a+b-1,0\}$. If $(X, {\cal F}, T_m)$ is a
complete Menger space and ${\cal C}$ is the collection of all
nonempty closed bounded subsets of $X$ in $(t,\epsilon)-$
topology, then $({\cal C}, {\cal F_C}, T_m)$ is also a complete
Menger space.
\end{prop}

{{\bf Proof.} See \cite{KS01}. $\hspace{1cm}\Box$ }

\subsection{ E-spaces}

\smallskip

\vspace{0.5cm} The notion of E-space was introduced by Sherwood \cite{Sh69}
in 1969. Next we recall this definition. Let $(\Omega ,{\cal K},P)$ be a
probability space and let $(Y,\rho )$ be a metric space. The ordered pair $(%
{\cal E,F})$ is an {\em E-space over the metric space $(Y,\rho )$} (briefly,
an E-space) if the elements of ${\cal E}$ are random variables from $\Omega $
into Y 
and ${\cal F}$ is the mapping from ${\cal E}\times {\cal E}$ into ${\Delta
^{+}}$ defined via ${\cal F}(x,y)=F_{x,y}$, where
\[
F_{x,y}(t)=P(\{\omega \in \Omega |\;d(x(\omega ),y(\omega ))<t\})
\]
for every $t\in $ ${\bf R}$. Usually $(\Omega ,{\cal K},P)$ is called the
base and $(Y,\rho )$ the target space of the E-space. If ${\cal F}$
satisfies the condition
\[
{\cal F}(x,y)\not=H,\;\;for\;\;x\not=y,
\]
with H defined in paragraf 3.1., then $({\cal E},{\cal F})$ is said to be a
{\em canonical E-space.} Sherwood \cite{Sh69} proved that every canonical
{\cal E}-space is a Menger space under $T=T_{m}$, where $T_{m}(a,b)=\max
\{a+b-1,0\}$. In the following we suppose that ${\cal E}$ is a canonical
E-space.

The convergence in an {\cal E}-space is exactly the probability convergence.
The E-space $({\cal E},{\cal F})$ is said to be complete if the Menger space
$({\cal E},{\cal F},T_{m})$ is complete.

\begin{prop}
\label{vulik}
If $ (Y,\rho)$ is a complete metric space then the  E-space $({\cal E,F}%
)$ is also complete.
\end{prop}

{\bf Proof.} This property is well-known for $Y=R$ (see e.g. \cite{V76},
Theorem VII.4.2.]). In the general case the proof is analogous.

Let $(x_n)_{n\in{\bf N}}$ be a Cauchy sequence of elements of ${\cal E}$,
that is
\[
\lim_{n,m\rightarrow \infty}F_{x_n,x_{n+m}}(t)=1,\,\mbox{for all} \, t>0.
\]

First we show that there exists a subsequence $(x_{n_k})_{k\in{\bf N}}$ of
the given sequence which is convergent almost everywhere to a random
variable x. Let as set positive numbers $\epsilon_i$ so that $%
\sum_{i=1}^\infty \epsilon_i<\infty$ and put $\delta_p=\sum_{i=p}^\infty%
\epsilon_i,\,p=1,2,....$ For each i there is a natural number $k_i$, such
that
\[
P(\{\omega\in\Omega \vert
\rho(x_k(\omega),x_l(\omega))\geq\epsilon_i\})<\epsilon_i\,\,\mbox{for}\,
k,l\geq k_i.
\]
We can assume that $k_1<k_2<...<k_i<...$. Then
\[
P(\{\omega\in\Omega \vert
\rho(x_{k_{i+1}}(\omega),x_{k_i}(\omega))\geq\epsilon_i\})<\epsilon_i\,\,%
\mbox{for}\, k,l\geq k_i.
\]
Let us put
\[
D_p=\cup_{i=p}^\infty\{\omega\in \Omega \vert\, \rho
(x_{k_{i+1}},x_{k_i})\geq \epsilon_i\}.
\]

Then $P(D_p)<\delta_p.$ Lastly, for the intersection $D^{\prime}=\cap_{p=1}^%
\infty D_p$ we obviously have $P(D^{\prime})=0$ since $\delta_p\rightarrow 0$%
. We shall show that the sequence $(x_{k_i}(\omega))$ has a finite limit $%
x(\omega)$ at every point ${\bf \omega}\in\{\omega\in\Omega \vert \,
\rho(x_k(\omega),x_m(\omega))>t\}\setminus D^{\prime}.$ For some p we have $%
x\notin D_p$. Consequently, $\rho(x_{k_{i+1}}(\omega),x_{k_i}(\omega))<%
\epsilon_i$, for all $i\geq p.$ It follows that for any two indices i and j
such that $j>i\geq p$ we have
\[
\rho(x_{k_{j}}(\omega),x_{k_i}(\omega))\leq
\sum_{m=i}^{j-1}\rho(x_{k_{m+1}}(\omega),x_{k_m}(\omega))<
\]
\[
< \sum_{m=i}^{j-1}\epsilon_m < \sum_{m=i}^{\infty}\epsilon_m=\delta_i .
\]
Thus $\lim_{i,j\rightarrow
\infty}\rho(x_{k_{j}}(\omega),x_{k_i}(\omega)))=0. $ This means that $%
(x_k(\omega))_{k\in{\bf N}}$ is a Chauchy sequence for every $\omega$ which
implies the pointwise convergence of $(x_{k_i})_{i\in{\bf N}}$ to a finite
limit function. Now it only remains to put
\[
x(\omega)=\left\{
\begin{array}{lcc}
\lim x_{k_i}(\omega) & for & \omega\notin D^{\prime} \\
0 & for & \omega\in D^{\prime}
\end{array}
\right.
\]
to obtain the desired limit random variable. By Lebeque theorem (see e.g.
\cite{V76} theorem VI.5.2) $x_{k_i}\rightarrow x$ with respect to d. Thus,
every Cauchy sequence in ${\cal E}$ has a limit, which means that the space $%
{\cal E}$ is complete. $\hspace{1cm}\Box $

\bigskip The next result was proved in \cite{KS01}:

\begin{th}
 Let $({\cal E, F})$ be a complete  E- space, $N\in {\bf N^*}$, and let $f_{1},...,f_{N}:
{\cal E}\rightarrow {\cal E}$ be  contractions with ratio
$r_{1},...r_{N}$, respectively. Suppose that there exists an
element $z\in {\cal E}$ and a real number $\gamma$ such that
\begin{equation}
\label{fotet}  P(\{\omega\in \Omega \vert
\rho(z(\omega),f_i(z(\omega))\geq t\})\leq\frac{\gamma}{t},
\end{equation}
 for
all $i\in\{1,..,N\}$ and for all $t>0.$  Then there exists a
unique nonempty closed bounded and compact subset $K$ of ${\cal
E}$ such that
\[
f_{1}(K)\cup...\cup f_{N}(K)=K.
\]
 \label{fo}
\end{th}

\begin{cor}\label{lem}
 Let $({\cal E, F})$ be a complete  E- space,  and let $f:{\cal
E}\rightarrow {\cal  E}$ be a contraction with ratio  $r$.
 Suppose there exists  $z\in {\cal E}$ and a real number $\gamma$ such that
  $$P(\{\omega\in\Omega\vert\; \rho(z(\omega),f(z)(\omega))\geq
  t\})\leq\frac{\gamma}{t}\;  \mbox{for all}\; t>0.$$
 Then there exists  a unique $x_0\in {\cal E}$  such that $f(x_0)=x_0.$
\end{cor}


\section{Proof of Theorem \ref{sajatmert}}

First we give two lemmas. Let ${\cal E}_q$ be the set of random variables
with values in $M_q$ and let ${\cal E}_q(\alpha)$ be the set
\[
{\cal E}_q(\alpha):=\{\beta\in {\cal E}_q\vert \, \exists \gamma>0\,
P(\{\omega\in \Omega\vert l_q(\alpha (\omega),\beta(\omega))\geq t\})\leq
\frac{\gamma}{t},\mbox{for all}\, t>0 \}.
\]

\begin{lem}
${\bf M}_q\subset {\cal E}_q(\alpha)$ for all $\alpha\in M_q.$
\end{lem}

{\bf Proof:} For $\beta\in {\bf M}_q$ we have
\[
P(\{\omega\in \Omega\vert l_q(\alpha(\omega),\beta(\omega))\geq
t\})=\int_{l_q(\alpha(\omega),\beta(\omega))\geq t}dP\leq
\]
\[
\leq \frac{1}{t}\int_\Omega l_q(\alpha(\omega),\beta(\omega))dP =\frac{1}{t}%
E_\omega l_q(\alpha(\omega),\beta(\omega)).
\]
Hence $\beta\in {\cal M}_q$ we have $\gamma=E_\omega
l_q(\alpha(\omega),\beta(\omega))<\infty$ for all $t>0.$ $\hskip1cm \Box$

\begin{lem}
$({\cal E}_q, {\cal F})$ is a complete E-space.
\end{lem}

{\bf Proof:} Choose $Y:={\cal E}_q$ and ${\cal F}_{\mu,\nu}
(t):=P(\{\omega\in \Omega\vert l_q(\mu(\omega),\nu(\omega))<t\})$ in the
Proposition \ref{vulik}. $\hskip1cm \Box$

{\bf Proof of Theorem \ref{sajatmert}:}
Let ${\cal S}$ be a random scaling law. Define $f:{\cal E}_q\rightarrow
{\cal E}_q$ by $f(\mu)={\bf S}\mu$, i.e.
\[
{\bf S}\mu(\omega)=\sum_ip_i^\omega S_i^\omega\mu(\omega^{(i)}).
\]
We first claim that if $\mu\in {\cal E}_q$ then ${\bf S}\mu\in {\cal E}_q.$
For this, choose iid $\mu(\omega^{(i)})\stackrel{d}{=}\mu(\omega)$ and $%
(p_1^\omega,S_1^\omega,...,p_N^\omega,S_N^\omega)\stackrel{d}{=}{\bf S }$
independent of $\mu(\omega).$ For $q\geq 1$ and $b=S_i^\omega(a)$ we compute
\[
\int d^q(x,a)d({\bf S}\mu^(\omega)(x))=l_q^q(\sum_{i=1}^Np_i^\omega
S_i^\omega\mu (\omega^{(i)}),\delta_a)=
\]
\[
= l_q^q(\sum_{i=1}^Np_i^\omega
S_i^\omega\mu(\omega^{(i)}),\sum_{i=1}^Np_i^\omega S_i^\omega\delta_b) \leq
\]
\[
\leq \sum_{i=1}^Np_i^\omega r_i^q l_q^q(\mu(\omega^{(i)}),\delta_b).
\]
Since $\mu (\omega^{(i)})\in M_q$ we have
\[
\int d^q(x,a)d({\bf S}\mu(x)< \infty.
\]

The case $0<q<1$ is dealt similarly, replacing $l_q^q$ by $l_q$:
\[
\int d^q(x,a)d({\bf S}\mu^(\omega)(x))=l_q(\sum_{i=1}^Np_i^\omega
S_i^\omega\mu (\omega^{(i)}),\delta_a)=
\]
\[
= l_q(\sum_{i=1}^Np_i^\omega
S_i^\omega\mu(\omega^{(i)}),\sum_{i=1}^Np_i^\omega S_i^\omega\delta_b) \leq
\]
\[
\leq \sum_{i=1}^Np_i^\omega r_i^q l_q(\mu(\omega^{(i)}),\delta_b)<\infty.
\]

To establish the contraction property let $\mu,\nu\in{\cal E}_q$, $%
\mu(\omega^{(i)})\stackrel{d}{=}\mu(\omega), \nu(\omega^{(i)})\stackrel{d}{=}%
\nu(\omega), i\in\{1,2,...,N\}$ and $q\geq 1$. We have
\[
F_{f(\mu),f(\nu)}(t)=P(\{\omega\in \overline \Omega\, \vert\,
l_q(f(\mu(\omega)),f(\nu(\omega))<t\})=
\]
\[
=P(\{\omega\in \overline \Omega \,\vert\, l_q(\sum_{i=1}^Np_i^\omega
S_i^\omega\mu(\omega^{(i)}), \sum_{i=1}^Np_i^\omega
S_i^\omega\nu(\omega^{(i)}))<t\})\geq
\]
\[
\geq P(\{\omega\in \overline \Omega \,\vert\, [\sum_{i=1}^Np_i^\omega
(r_i)^q l_q^q(\mu(\omega^{(i)}),\nu(\omega^{(i)}))]^\frac 1q<t\})\geq
\]
\[
\geq P(\{\omega\in \overline \Omega \,\vert\, [\lambda_q
l_q^q(\mu(\omega),\nu(\omega))]^\frac 1q<t\})=F_{\mu,\nu}(\frac{t}{%
\lambda_q^{\frac 1q}})
\]
for all $t>0$. In case $0<q<1$, one replaces $l_q^q$ everywhere by $l_q$.
Thus ${\bf S }$ is a contraction map with ratio $\lambda_q^{\frac{1}{q}%
\wedge 1}$. We can apply Corollary \ref{lem} for $r=\lambda_q^{\frac
1q\wedge 1}$. If $\mu^*$ is the unique fixed point of ${\bf S}$ and $%
\mu_0\in M_q$ 
then
\[
F_{{\bf S}^n\mu_0,\mu^*}(t)=P(\{\omega\in\overline\Omega\,\vert\, l_q({\bf S}%
^n\mu_0,\mu^*)<t\})\geq
\]
\[
\geq P(\{\omega\in\overline\Omega\,\vert\, \frac{\lambda_q^{\frac nq}}{%
1-\lambda_q^{\frac 1q}}l_q(\mu_0,{\bf S}\mu_0)<t\}).
\]
and
\[
\lim_{n\rightarrow \infty} F_{{\bf S}^n\mu_0,\mu^*}(t)=1\,\mbox{for all }\,
t>0.
\]

From $\mu_{n+1}(\omega)={\bf S}\mu_n(\omega)$ it follows that $%
\mu_m\rightarrow \mu^*$exponentially fast. Moreover, for $q\geq 1$
\[
\sum_{i=1}^\infty \overline P(l_q^q({\bf S}^n\nu_0,\mu^*)\geq\epsilon)\leq
\sum_{i=1}^\infty \frac{el_q^q({\bf S}^n\mu_0,\mu^*)}{\epsilon}\leq
c\sum_{i=1}^\infty \frac{\lambda_q^n}{\epsilon}<\infty.
\]
This implies by Borel Catelli lemma that $l_q(\mu_n,\mu^*)\rightarrow 0$ a.s.

For the uniqueness let ${\cal D}$ the set of probability distribution of
members of ${\cal E}_q.$ We define on ${\cal D}$ the probability metric by
\[
F_{{\cal A,B}}(t)=\sup_{s<t}\sup\{F_{\mu,\nu}(s)\vert \,\mu\stackrel{d}{=}%
{\cal A},\nu\stackrel{d}{=}{\cal B}\}.
\]
To establish the contraction property, let ${\cal A,B\in D}$. For $q\geq 1$,
on has
\[
F_{{\cal SA,SB}}(t)=\sup_{s<t}\sup\{F_{{\bf S}\mu,{\bf S}\nu}(s)\vert\, \mu%
\stackrel{d}{=}{\cal A},\nu\stackrel{d}{=}{\cal B}\}\geq
\]
\[
\geq\sup_{s<t}\sup\{F_{\mu,\nu}(\frac{s}{\lambda_q})\vert \, \mu\stackrel{d}{%
=}{\cal A},\nu\stackrel{d}{=}{\cal B}\}= F_{{\cal A,B}}(\frac{t}{\lambda_q})
\]
for all $t>0.$ In case $0<q<1$ on work similarly.

Let ${\cal D}_1$ and ${\cal D}_2$ such that ${\cal SD}_1={\cal D}_1$ and $%
{\cal SD}_2={\cal D}_2$. 

Since ${\cal D}_1={\cal S}^n({\cal D}_1)$ and ${\cal D}_2={\cal S}^n({\cal D}%
_2)$ we have
\[
F_{{\cal D}_1,{\cal D}_2}(t)\geq F_{{\cal D}_1,{\cal D}_2}(\frac {t}{r^n})
\]
for all $t>0$. Using $\lim_{n\rightarrow \infty}r^n=0$ it follows that

\[
F_{{\cal D}_{1},{\cal D}_{2}}(t)=1,
\]
for all $t>0.$ $\hspace{1cm}\Box $

{\bf Remark:} Since $\lambda_q^{\frac 1q}\rightarrow max_ir_i$ as $%
q\rightarrow \infty$, we can regard Theorem 3.1. from \cite{KS01} as a limit
case of Theorem \ref{sajatmert}. More precisely, if $max_ir_i<1$ then $sprt
\mu^*$ is the unique compact set satisfying $(S_1,...,S_N).$

\bigskip

\vskip2cm

\end{document}